\documentclass[11pt,numreferences,draft]{kluwer}
\usepackage{mathrsfs}
\usepackage{amssymb}

\usepackage[english]{babel}
\newtheorem{theo}{Theorem}[section]
\newtheorem{lema}{Lemma}[section]
\newtheorem{prop}{Proposition}[section]

\newcommand{\bq}{\begin{equation}}
\newcommand{\eq}{\end{equation}}
\newcommand{\ba}{\begin{array}}
\newcommand{\ea}{\end{array}}

\newcommand{\pe}[2]{\langle#1,#2\rangle}

\newcommand\bd{\begin{definicion}{\bf }}
\newcommand\ed{\end{definicion}}
\newcommand\bl{\begin{lema}{\bf }}
\newcommand\el{\end{lema}}
\newcommand\bp{\begin{prop}{\bf }}
\newcommand\ep{\end{prop}}
\newcommand\bt{\begin{theo}{\bf }}
\newcommand\et{\end{theo}}
\newcommand\bdm{\begin{proof}}
\newcommand\edm{\end{proof}}
\newcommand\bn{\begin{nota}{\bf }}
\newcommand\en{\end{nota}}
\newcommand\bc{\begin{corolary}{\bf }}
\newcommand\ec{\end{corolary}}

\newtheorem{nota}{Remark}[section]

\newtheorem{definicion}{Definition}[section]
\newtheorem{corolary}{Corollary}[section]
\newenvironment{proof}{{\bf Proof:\ }}{\hfill
$\Box$}
\newcommand\eqref[1]{(\ref{#1})}

\newfont{\got}{eufm10 scaled \magstep1}

\date{\today}
\begin{document}
\begin{article}
\begin{opening}
\title{$q$-Classical orthogonal polynomials:
A general difference calculus approach}
\author{R. S. \surname{Costas-Santos}
\email{rscosa@gmail.com}}
\institute{Department of Mathematics. University of
California, Santa Barbara, CA 93106, USA.}
\author{F. \surname{Marcell\'{a}n}
\email{pacomarc@ing.uc3m.es}}
\institute{Departamento de Matem\'{a}ticas. Universidad Carlos III de Madrid.
Avenida Universidad 30, 28911-Legan\'{e}s, Spain}
\date{\today}
\runningauthor{R. S. Costas-Santos and F. Marcell\'{a}n}
\noindent 
\begin{abstract}
It is well known that the classical families
of orthogonal polynomials are characterized
as the polynomial eigenfunctions of a second
order homogenous linear differential/difference
hypergeometric operator with polynomial coefficients.
\\
In this paper we present a study of the classical
orthogonal polynomials sequences, in short
classical OPS, in a more general framework by
using the differential (or difference) calculus
and Operator Theory.
\\
The Hahn's Theorem and a characterization
Theorem for the $q$-polynomials which belongs
to the $q$-Askey and Hahn tableaux are proved.
Finally, we illustrate our results applying them
to some known families of orthogonal $q$-polynomials.
\end{abstract}
\classification{2000 MSC codes}{33C45; 33D45 }
\end{opening}
\section{Introduction} \label{sec1}
Classical orthogonal polynomials constitute
a very important and interesting family of
special functions.
They are mathematical objects which
have attracted attention not only because
of their mathematical value but also
because of their connections with
physical problems.
In fact, they are also related, among
others, to continued fractions, Eulerian
series, elliptic functions \cite{and,fin},
and quantum algebras \cite{koo1,koo2,vikl}.
\\
They also satisfy a three-term recurrence
relation (TTRR) \cite{nisuuv}
$$
x(s)p_n(x(s))=\alpha_n p_{n+1}(x(s))+
\beta_n p_n(x(s))+\gamma_n p_{n-1}(x(s)),
\quad n\ge 0,
$$
where $\gamma_n\ne 0$, $n\ge 1$, as well as
their derivatives (or differences or
$q$-differences) also constitute a a sequence
of orthogonal family (see e.g.
\cite{alv2,alv3,atrasu,nisuuv} for a more
recent review).

Indeed, a fundamental role is played by the
so-called {\sf characterization Theorems},
i.e. the Theorems which collect those properties
that completely define and characterize the
classical orthogonal polynomials.
\\
One of the many ways to characterize a family
$(p_n)$ of classical polynomials
(Hermite, Laguerre, Jacobi, and Bessel), which
was first posed by R. Askey and proved by
W. A. Al-Salam and T. S. Chihara \cite{alch}
(see also \cite{mabrpe}), is the structure
relation
\begin{equation}
\phi(x)p'_n(x)=\tilde a_np_{n+1}(x)+\tilde b_n
p_n(x)+\tilde  c_np_{n-1}(x), \quad n\ge 0,
\label{1:1}
\end{equation}
where $\phi$ is a fixed polynomial of degree at
most 2 and $\tilde c_n\ne 0$, $n\ge 1$.
\\
A. G. Garcia, F. Marcell\'{a}n, and L. Salto
\cite{gamasa} proved that the relation
\eqref{1:1} also characterizes the discrete
classical orthogonal polynomials (Hahn,
Krawtchouk, Meixner, and Charlier polynomials)
when the  derivative is replaced by the forward
difference operator $\Delta$ defined as
$$
\Delta f(x):=({\mathscr E}^+-\mathscr{I})f(x)=
f(x+1)-f(x),
$$
where ${\mathscr E}^+$ is the forward shift operator
and ${\mathscr I}$ the identity operator.

Later on, J. C. Medem, R. \'{A}lvarez-Nodarse,
and F. Marcell\'{a}n \cite{mealma} (see also
\cite{alal,alv3,alar1}) characterized the
orthogonal polynomials which belong to the
$q$-Hahn class by a structure relation obtained
from \eqref{1:1} replacing the derivative by the
$q$-difference operator ${\mathscr D}_q$ defined
as
$$
{\mathscr D}_q(f)(x):=\frac{f(qx)-f(x)}{(q-1)x},
\quad q \in {\mathbb C}, |q|\ne 0,1.
$$
The aim of the present paper is to go further in
the study of classical polynomials giving two
new structure relations for the $q$-polynomials
which belong to the $q$-Askey tableau and
the lattice is $q$-quadratic, i.e. $x(s)=c_1q^s+
c_2q^{-s}+c_3$, with $c_1c_2\ne 0$, being
$q\in \widetilde {{\mathbb C}}:={\mathbb C}
\setminus\left(\{0\}\cup\left(\bigcup_{n\ge 1}\{
z\in {\mathbb C}: z^n=1\}\right)\right).$

In fact, we prove that the following relation
characterizes the $q$-polynomials.
$$
{\mathscr M} p_{n}(x(s))=
e_n\frac{\Delta p_{n+1}(s)} {\Delta x(s)}+
f_n \frac{\Delta p_{n}(s)}{\Delta x(s)}
+g_n\frac{\Delta p_{n-1}(s)}{\Delta x(s)},
\quad n\ge 0,
$$
where ${\mathscr M}$ is the forward arithmetic
mean operator:
$$
{\mathscr M}:=\frac{1}{2}({\mathscr E}^++\mathscr{I}),
$$
and, $(e_n)$, $(f_n)$, and $(g_n)$ are sequences of complex numbers
such that $e_n\ne 0$ and $g_n\ne \gamma_n$.
\\
Furthermore, we give a characterization
Theorem for classical orthogonal polynomials
in a more general framework using Operator
Theory, as well as we deduce some unified
expressions for the second order homogeneous
linear differential (or difference)
hypergeometric operator and
for its polynomial eigenfunctions.
\\
The structure of this paper is as follows: In
Section 2 we introduce some notations and
definitions used throughout the paper.
In Section 3 we define the Rodrigues Operator
and using Operator Theory we deduce a unified
expression for the linear
differential (difference resp.) hypergeometric
operators related to the classical families, and
for their polynomial eigenfunctions, as well as
other algebraic properties.
In Section 4 we present a new characterization
Theorem for $q$-classical orthogonal polynomials
and we prove the Hahn's Theorem for the
$q$-polynomials of the $q$-Askey tableau.
In Section 5 we illustrate our results applying
them to some known $q$-polynomials.
\section{Preliminaries} \label{sec2}
The standard classical orthogonal polynomials
(Hermite, Laguerre, Jacobi, and Bessel) are
eigenfunctions of a second order linear
differential operator \cite{routh}
\begin{equation} \label{2:1}
{\mathscr H}:=\tilde{\sigma}(x)\frac{d^2}{dx^2}+
\tilde{\tau}(x)\frac{d}{dx},
\end{equation}
where $\tilde{\sigma}$ and $\tilde{\tau}$ are
polynomials of degree at most 2 and 1,
respectively.
This operator is also said to be a {\sf standard
Hamiltonian operator} (see \cite{nisuuv}).
\\
In fact, a sequence of monic polynomials,
$(p_n)$, such that
$$
{\mathscr H} p_n(x)=\left(n\tilde{\tau}'+n(n-1)
\frac{\tilde{\sigma}''}2\right)p_n(x), \quad n\ge 0,
$$
is orthogonal with respect to a certain weight
function $\rho(x)$ supported on $\Omega\subseteq
{\mathbb R}$, i.e.
\begin{equation}
\displaystyle \int_\Omega p_n(x)p_m(x)\rho(x)
dx=d_n^2 \delta_{n,m}, \quad n, m\ge 0. \label{2:3}
\end{equation}
Notice that
$\Omega$ is an interval associated with the regular
singularities of ${\mathscr H}$.

Moreover, this weight function $\rho(x)$ fulfills
the following Pearson equation
$$
\frac{d}{dx}(\tilde{\sigma}(x)\rho(x))=
\tilde{\tau}(x)\rho(x).
$$
On the other hand, we can consider two
different discretizations of the operator
\eqref{2:1} replacing the derivative by
certain approximations on a lattice. Thus
we need first to introduce the concept of lattice.
\bd \label{def2.1}\cite{atrasu,nisuuv}
A {\sf lattice} is a complex function $x\in
C^2(\Lambda)$ where $\Lambda$ is a complex domain,
${\mathbb N}_0\subseteq \Lambda$, and $x(s)$, $s=0, 1,
\dots$ are the points where we will discretize
\eqref{2:1}.
\ed
For a first approximation of the operator
${\mathscr H}$, we consider the uniform lattice
$x(s)=s$.
Thus the following second order homogeneous linear
difference hypergeometric operator
\begin{equation} \label{2:4}
{\mathscr H}_{\Delta}:=\displaystyle \sigma(s)
\Delta\nabla+\tau(s)\Delta,
\end{equation}
can be introduced where $\nabla$ and $\Delta$ are
the backward and forward difference operators,
respectively, being
$$
\nabla f(x)=({\mathscr E}^--\mathscr{I})f(x)
=f(x)-f(x-1),
$$
$\sigma(x):=\tilde{\sigma}(x)-\frac{1}{2}
\tilde{\tau}(x)$ is a polynomial of degree at
most 2, and $\tau(x)=\tilde{\tau}(x)$.
\\
Indeed, a sequence of polynomials $(p_n)$
satisfying
$$
{\mathscr H}_\Delta p_n(x)=\left(n\tilde{\tau}'+
n(n-1)\frac{\tilde{\sigma}''}2\right) p_n(x),
\quad n\ge 0,
$$
is orthogonal with respect to a certain weight
function $\rho(x)$ supported on $\Omega\subseteq
{\mathbb R}$, i.e.
\begin{equation} \label{2:5}
\sum_{x\in \Omega} p_n(x)p_m(x)\rho(x)=
d_n^2 \delta_{n,m}, \quad \quad n, m\ge 0,
\end{equation}
with some extra boundary conditions (see
\S 2.3 in page 26 of \cite{nisuuv}).

As in the continuous case, the weight function
$\rho$ fulfills the following $\Delta$-Pearson
equation
$$
\Delta(\sigma(x)\rho(x))=\tau(x)\rho(x).
$$
The second way to discretize the operator
${\mathscr H}$ is to consider a nonuniform
lattice $x(s)$, (see \cite{nisuuv}) which
leads to the following linear difference hypergeometric
operator:
\begin{equation} \label{2:6}
{\mathscr H}_q:=\sigma(s)\frac{\Delta}{\nabla
x_1(s)}\frac{\nabla} {\nabla x(s)}+\tau(s)
\frac{\Delta}{\Delta x(s)},
\end{equation}
where $x(s)=c_1q^s+c_2q^{-s}+c_3$,
$q\in \widetilde {{\mathbb C}}$,
$\sigma(s):=\tilde{\sigma}(x(s))-\frac{1}{2}
\tau(x(s))\nabla x_1(s)$, $x_1(s):=x(s+\frac{1}
{2})$, and $\tau(s)=\tilde{\tau}(x(s))$.

Notice that, in general, $\sigma(s)$ is not a
polynomial on $x(s)$ but if either
$c_1=0$ or $c_2=0$ then $\sigma(s)$ is
a polynomial on $x(s)$ of degree at most 2,
i.e. $\sigma(s)\equiv\sigma(x(s))$.

Indeed, a sequence of monic polynomials
$(p_n)$ satisfying
$$
{\mathscr H}_q p_n(x(s))=[n]_q\left(\alpha_q
(n-1)\,\tilde{\tau}'+[n-1]_q
\frac{\tilde{\sigma}''}2 \right) p_n(x(s)),
\quad n\ge 0,
$$
where $[s]_q$ denotes the $q$-number
(in its symmetric form)
$$
[s]_q:=\frac{q^{\frac s2}-q^{-\frac s2}}
{q^{\frac 12}-q^{-\frac 12}},\qquad s\in
{\mathbb C},
$$
and
$$
\alpha_q(s):=\frac{q^{\frac {s}2}+
q^{-\frac{s}2}}2, \quad s\in {\mathbb C},
$$
is orthogonal with respect to a certain
linear functional $\bf u\in {\mathbb P}'$, i.e.
$$
\pe {{\bf u}} {p_np_m}=d^2_n\delta_{m,n},
\quad m, n\ge 0.
$$
\bn \label{rem2.1}
Notice that $\bf u$ admits a Lebesgue-Stieltjes
integral representation, not necessarily unique,
$$
\pe {{\bf u}}p=\int_{\Gamma}p(z)d\mu(z),
$$
where $\Gamma$ is a contour in the complex
plane and $\mu$ is a measure, with
supp$(\mu)=\Gamma$.

In fact,  in most of the cases, such an integral
can be transformed into a sum in such a way
the linear functional can be written as \cite{kost}
\begin{equation} \label{2:10}
\pe{{\bf u}}{p}=\sum_{s_i\in \Omega} p(x(s_i))
\rho(s_i)\nabla x_1(s_i), \quad s_{i+1}=s_i+1,
\end{equation}
where the function $\rho$ is supported on
$\Omega\subseteq {\mathbb R}$.
For this reason we will assume the above
representation for the functional although
the proofs when $\bf u$ can not be
represented in such a way are analogous and we will omit them.
\en

Moreover, as in the discrete case, the weight
function $\rho(s)$ satisfies the following
$q$-Pearson equation
\begin{equation} \label{2:11}
\frac{\Delta\big(\sigma(s)\rho(s)\big)}
{\nabla x_1(s)}=\tau(s)\rho(s).
\end{equation}
\bn \label{rem2.2}
It is important to point out that in the
orthogonality conditions \eqref{2:3},
\eqref{2:5}, and \eqref{2:10} one
needs to impose some boundary conditions,
which follow from the fact that the
moments of the corresponding measure must
be finite (see, for instance, \cite{ch},
\cite[p. 27, p. 70]{nisuuv},
\cite{niuv1}).
\en
An important fact related to these operators
is its hypergeometric character, i.e. if
$y$ is an eigenfunction of the linear
operator ${\mathscr H}$, ${\mathscr H}_\Delta$,
or ${\mathscr H}_q$ respectively,
its derivative $\tilde y=y'$, its difference
$\tilde y=\Delta y$,  or its $q$-difference
$\tilde y=\Delta^{(1)}y$, respectively,
is an eigenfunction of the second order linear
operator of the same kind.
For instance, $\tilde y=\Delta^{(1)}y$ is an
eigenfunction of the following  second order
linear difference operator:
$$
\widetilde{{\mathscr H}_q}:=\sigma(s)\frac
{\Delta}{\nabla x_1(s)}\frac{\nabla}{\nabla
x(s)}+\left(\tau(s+1)\frac{\Delta x_1(s)}
{\Delta x(s)}+\frac{\sigma(s+1)-\sigma(s)}
{\Delta x(s)}\right)\frac{\Delta}{\Delta x(s)},
$$
where
$$
\Delta^{(m)}:=\frac{\Delta}{\Delta x_{m-1}(s)}
\frac{\Delta}{\Delta x_{m-2}(s)}\dots
\frac{\Delta}{\Delta x(s)}, \quad m\ge 1,
\quad x_m(s):=x(s+\mbox{$\frac m 2$}), \quad
m\in{\mathbb Z}.
$$
Throughout the paper we will study some
properties related to the difference
operator ${\mathscr H}_q$ and its eigenfunctions.

We can rewrite ${\mathscr H}_q$ as follows
\begin{equation} \label{2:12}
{\mathscr H}_q=\frac 1 {\nabla x_1(s)}
\left(\Big(\sigma(s)+\tau(s)\nabla x_1(s)
\Big) \frac{\Delta}{\Delta x(s)}-\sigma(s)
\frac{\nabla}{\nabla x(s)}\right).
\end{equation}
Moreover, with this definition we can also
rewrite the $q$-Pearson equation \eqref{2:11}
as
\begin{equation} \label{2:13}
\sigma(s+1)\rho(s+1)=\Big(\sigma(s)+\tau(s)
\nabla x_1(s)\Big)\rho(s).
\end{equation}
Combining \eqref{2:12} and
\eqref{2:13}, we get the symmetric or
self-adjoint form of \eqref{2:6}
\begin{equation} \label{2:14}
{\mathscr H}_q=\left[\frac 1 {\rho(s)}\frac{\nabla}
{\nabla x_1(s)}\, \rho_1(s)\right]\circ \frac{\Delta}
{\Delta x(s)},
\end{equation}
where $\rho_0(s):=\rho(s)$ and $\rho_1(s):=
\rho(s+1)\sigma(s+1)$.
\\
This representation for the difference operator
\eqref{2:6} is the key to define the
Rodrigues Operator which we will use
in order to unify all the representations in
the context of the classical orthogonal
polynomials.
\bn \label{rem2.3}
A direct calculation yields the analogous formulas
for the operators ${\mathscr H}$ and
${\mathscr H}_{\Delta}$.
$$
{\mathscr H}=\left[\frac 1 {\rho(x)} \frac d{dx}
\rho_1(x)\right] \circ \frac{d} {dx},
$$
where $\rho_0(x):=\rho(x)$ and $\rho_1(x):=
\rho(x)\sigma(x)$; and
$$
{\mathscr H}_\Delta=\left[\frac 1 {\rho(s)}
\nabla \rho_1(s)\right]\circ \Delta,
$$
where $\rho_0(s):=\rho(s)$ and $\rho_1(s):=
\rho(s+1)\sigma(s+1)$.
\en
Using this information we can
write the polynomial eigenfunctions
of the operator ${\mathscr H}_q$ as
(see \cite[(3.2.10), p. 66]{nisuuv})
\begin{equation} \label{2:16}
p_n(x(s)):=\left[\frac {B_n}{\rho_0(s)}
\frac{\nabla} {\nabla x_1(s)}\rho_1(s)
\right]\circ \left[\frac 1{\rho_1(s)}
\frac{\nabla}{\nabla x_2(s)}\rho_2(s)
\right] \circ  \cdots \circ \left[\frac
1{\rho_{n-1}(s)}\frac{\nabla}
{\nabla x_n(s)}\, \rho_n(s)\right],
\end{equation}
or, after some simplifications, as
$$
p_n(x(s))=\frac {B_n}{\rho_0(s)}
\nabla_0^{(n)}\, [\rho_n(s)], \quad
B_n \ne 0,
$$
where $\rho_0(s)=\rho(s)$, $\rho_k(s)=
\rho_{k-1}(s+1)\sigma(s+1)$, for every
$k\in \mathbb Z$, and
$$
\nabla^{(n)}\equiv \nabla_{0}^{(n)},\quad
\nabla_{k}^{(n)}:=\frac{\nabla} {\nabla
x_{k+1}(s)}\frac{\nabla}{\nabla x_{k+1}(s)}
\cdots \frac{\nabla}{\nabla x_n(s)},\qquad
k=0, 1, \dots, n-1.
$$
\subsection{The (forward) mean arithmetic process}
\label{sec2.1}
It is well known that for
any pair of polynomials, $\pi$ and $\xi$,
the following relations hold:
\begin{eqnarray}
\nonumber\Delta^{(1)}(\pi \xi)=&
(\Delta^{(1)}\pi) \xi+({\mathscr E}^+
\pi) \Delta^{(1)}\xi,\\ \Delta^{(1)}
(\pi \xi)=& (\Delta^{(1)}\xi) \pi+
({\mathscr E}^+\xi) \Delta^{(1)}\pi.
\nonumber
\end{eqnarray}
We will say that we apply the {\sf (forward)
mean arithmetic process to} $\Delta^{(1)}(\pi
\xi)$ if we get the arithmetic mean of the
above two expressions obtaining in such a case
$$
\Delta^{(1)}(\pi \xi)=(\Delta^{(1)} \pi)
{\mathscr M}\xi+(\Delta^{(1)}\xi){\mathscr
M}\pi.
$$
For example,
$$
\Delta^{(1)}(x(s)p_n(x(s)))={\mathscr M}
p_n(x(s))+\alpha_q(1)x_1(s)\Delta^{(1)}
p_n(x(s))+c_3(q)(1-\alpha_q(1)).
$$
In an analogous way we define the
(backward) mean arithmetic mean process
although we will not use it.
\section{The Rodrigues Operator} \label{sec3}
The expressions given in the previous
section for the difference operator
${\mathscr H}_q$ and for its polynomial
eigenfunctions \eqref{2:16} suggest us
to consider a new operator which we call
the {\sf Rodrigues Operator}.
\bd \label{def3.1}
Given functions $\sigma$ and $\rho$, where
$\rho$ is supported on $\Omega$, and a
lattice $x(s)$, we define the {\sf $k$-th
Rodrigues Operator associated with}
($\sigma(s)$, $\rho_m(s)$, $x_m(s)$)
as
$$ \ba{c}
\displaystyle {\mathscr R}_0(\sigma(s),
\rho_m(s), x_m(s)):={\mathscr I}, \qquad
{\mathscr R}_1(\sigma(s), \rho_m(s),x_m
(s)):=\frac 1 {\rho_m(s)}\frac{\nabla}
{\nabla x_{m+1}(s)}\, \rho_{m+1}(s)\circ
{\mathscr I}, \\[5mm] {\mathscr R}_k(
\sigma(s),\rho_m(s),x_m(s)):={\mathscr
R}_1(\sigma(s), \rho_m(s),x_m(s)) \circ
{\mathscr R}_{k-1}(\sigma(s),\rho_{m+1}
(s), x_{m+1}(s)), \ k\ge 2.
\ea $$
\ed
For the  sake of simplicity, we denote
$$
{\mathscr R}_k(\rho_m(s),x_m(s))
:={\mathscr R}_k(\sigma(s),\rho_m(s),
x_m(s)),\qquad m\in {\mathbb Z}, \ k\ge 0.
$$
Notice that in the continuous
case
$$
{\mathscr R}_1(\rho_m(x),x)=\frac 1{\rho_m(x)}
\frac d{dx}\rho_{m+1}(x)\circ {\mathscr I},
$$
and in the  discrete case
$$
{\mathscr R}_1(\rho_m(s),x(s))=\frac 1{\rho_m
(s)}\nabla\rho_{m+1}(s)\circ {\mathscr I}.
$$
From this definition the difference
operator \eqref{2:14} can be written as
\bq \label{aux1}
{\mathscr H}_q={\mathscr R}_1(\rho(s),
x(s))\circ \Delta^{(1)},
\eq
and \eqref{2:16} reads
$$
p_n(x(s))=B_n\, {\mathscr R}_n(\rho(s),
x(s))(1), \quad B_n\ne 0,\ n\ge 0.
$$
Therefore, this is the way to write the
operators ${\mathscr H}$ (${\mathscr
H}_\Delta$ or ${\mathscr H}_q$, resp.)
and their corresponding polynomial
eigenfunctions in a unified way.
\\
Now, let us consider other properties
related to the Rodrigues Operator.
\bl \label{lem3.1}
The Pearson-type equation \eqref{2:11}
can be written as follows
$$
{\mathscr R}_1(\sigma(s),\rho(s),
x(s))(1)=\tau(s).
$$
\el
As a consequence we get
result:
\bl \label{lem3.2}
\cite[(3.2.18), p.66]{nisuuv}
For every integers $n,\, k$, $0\le
k\le n$, there exists a constant,
$C_{n,k}$, such that
$$
\Delta^{(k)}p_n(x(s))=C_{n,k}
{\mathscr R}_{n-k} (\rho_{k}(s),x_k(s))(1).
$$
\el
\section{The characterization Theorem for
the $q$-classical polynomials} \label{sec4}
The characterization Theorems constitute a
very useful tool to analyze classical
orthogonal polynomials as we stated in
the introduction.
Recently, in \cite{alv3} and \cite{alal}
a characterization of $q$-classical
polynomials in the $q$-linear lattice
$x(s)=cq^{\pm s}+d$ is presented.
We will extend these results for the
$q$-quadratic lattice
\begin{equation} \label{4:1}
x(s)=c_1q^s+c_2q^{-s}+c_3,\qquad q\in
\widetilde{{\mathbb C}}.
\end{equation}
For our purpose, we need to redefine the
concept of $q$-classical OPS. Moreover,
because the measures are supported on
the real line, we will consider
$\Omega$ as the interval\footnote{If
$|b|<|a|=\infty$ (resp. $|a|<|b|=\infty$)
then we write $(a,b]$ (resp. $[a,b)$);
if $|b|=|a|=\infty$
then we write $(a,b)$.} $[a,b]$.
\bd \label{def4.1}
The sequence $(p_n)$, where $p_n$ is a
polynomial on $x(s)$ of degree $n$ for all
$n\in {\mathbb N}_0$, is said to be a
$q$-{\sf classical OPS} on the lattice
$x(s)$ if it satisfies the property
of orthogonality \eqref{2:10}
where
\begin{enumerate}
\item[(i)] $\rho(s)$ is a solution of
the $q$-Pearson equation
\begin{equation} \label{4:2}
\Delta(\sigma(s)\rho(s))=\tau(s)\rho(s)
\nabla x_1(s).
\end{equation}
\item[(ii)] $\sigma(s)+\frac{1}{2} \tau(s)
\nabla x_1(s)$ is a polynomial on $x(s)$
of degree, at most, 2.
\item[(iii)] $\tau(s)$ is a polynomial on
$x(s)$ of degree 1.
\end{enumerate}
$(\sigma(s),\rho(s),x(s))$ is said to be
{\sf $q$-classical} if it satisfies (i)-(iii).
\ed
Notice that the above definition coincides
with Definition 1.1 given in \cite{alv3}
for the $q$-linear lattices and generalizes it
to the $q$-quadratic lattice.
\bl \label{lem4.1} \cite[(3.2.5), p. 62]{nisuuv}
For every polynomial $\pi$ and every integer $k$,
$\pi(x_k(s))+\pi(x_k(s-1))$ is a polynomial on
$x_{k-1}(s)$ of degree $\deg(\pi)$.
\el
\bl \label{lem4.2}
If $(\sigma(s),\rho(s),x(s))$ is $q$-classical,
then for every integer $k$ and every polynomial
$\pi$, the function
$$
{\mathscr R}_1(\rho_{k}(s),x_k(s))
\big(\pi(x_{k+1}(s))\big)
$$
is a polynomial on $x_k(s)$ of degree $\deg(\pi)+1$.
\el
\bdm
By hypothesis, we know that
$\widehat{\sigma}(x(s)):=\sigma(s)+
\frac{1}{2}\tau(s)\nabla x_1(s)$ and
$\tau(s)$ are polynomials on $x(s)$
of degree, at most 2, and 1, respectively.
From a direct calculation we get that
$$
{\mathscr R}_1(\rho_k(s),x_k(s))(1)=
\frac{\sigma(s+k)+\tau(s+k)\nabla
x_1(s+k)-\sigma(s)}{\nabla x_{k+1}(s)},
$$
is a polynomial on $x_k(s)$ of degree 1.

Furthermore, for every integer $k$ and
every polynomial $\pi$, we get
$$ \ba{rl}
\displaystyle {\mathscr R}_1(\rho_k(s),
x_k(s)) \big(&\!\!\!\pi(x_{k+1}(s))
\big)=\displaystyle \frac{\big(\sigma(s+k)
+\tau(s+k)\nabla x_1(s+k)\big)} {\nabla
x_{k+1}(s)} \pi(x_{k+1}(s))-\frac{\sigma(s)}
{\nabla x_{k+1}(s)}\\[4mm] & \displaystyle
\times \pi(x_{k+1}(s-1)) =\frac{\Theta(s+k)}
{\nabla x_{k+1}(s)} \pi(x_{k+1}(s))-\frac{
\sigma(s)}{\nabla x_{k+1}(s)}\pi(x_{k+1}(s-1)),
\ea$$
being $\Theta(s):=\sigma(s)+ \tau(s)\nabla x_1(s)$.

So, ${\mathscr R}_1(\rho_k(s),x_k(s))
\big(\pi(x_{k+1}(s))\big)$ is a polynomial
on $x_k(s)$ if and only if
\begin{equation} \label{4:3}
\frac{\Theta(s+\mbox{$\frac k2$})}{\nabla
x_1(s)}\pi(x_1(s))-\frac{\sigma(s-
\mbox{$\frac k2$})}{\nabla x_1(s)}\pi(x_{-1}(s))
\end{equation}
is a polynomial on $x(s)$.
But we can rewrite this expression in
the form
$$
\frac{\Theta(s+\mbox{$\frac k2$})-
\sigma(s-\mbox{$\frac k2$})}
{\nabla x_1(s)}\pi(x_1(s)) +\frac{\sigma
(s-\mbox{$\frac k2$})}{\nabla x_1(s)}
\Big(\pi(x_1(s))-\pi(x_{-1}(s))\Big),
$$
or, equivalently,
$$
\frac{\Theta(s+\mbox{$\frac k2$})-
\sigma(s-\mbox{$\frac k2$})}{\nabla
x_1(s)}\pi(x_{-1}(s))+\frac{\Theta(s+
\mbox{$\frac k2$})}{\nabla x_1(s)}
\Big(\pi(x_1(s))-\pi(x_{-1}(s))\Big).
$$
Taking the arithmetic mean of the above
expressions, using Lemma \ref{lem4.1}
as well as the relation
$$
\Theta(s+\mbox{$\frac k2$})-\sigma(s-
\mbox{$\frac k2$})=q_1(x(s))\nabla x_1(s),
$$
where $q_1$ is a polynomial of degree 1,
we deduce that \eqref{4:3} is a polynomial
on $x(s)$ of degree, at most, $\deg(\pi)+1$.
Moreover, from a straightforward calculation
the coefficient of $x^{m+1}$ is
$$
\alpha_q(m+2k)\, \tau'+[m+2k]_q
\frac{\widehat \sigma''}{2}\ne 0,
\quad m=\deg(\pi),
$$
where $\widehat \sigma(x)=\frac{\widehat
\sigma''}{2}x^2+\widehat \sigma'(0)x+
\widehat \sigma(0)$ and $\tau(x)=
\tau'x+\tau(0)$.
Hence the result follows.
\edm
\\
Taking into account all these results, we
can  state the preliminary
results related to our main Theorem.
\bt \label{theo4.1}
Let $(p_n)$ be a $q$-classical
OPS with respect to $\rho(s)$ on $x(s)$
such that
\begin{equation} \label{4:4}
x^k(a)x^l_{-1}(a)\sigma(a)\rho(a)=
x^k(b)x^l_{-1}(b)\sigma(b)\rho(b)=0,
\quad  k,l\ge 0.
\end{equation}
Then, $(\Delta^{(1)}p_{n+1})$
is a $q$-classical OPS with respect to the
function $\rho_1(s)$ on $x_1(s)$.
\\
Furthermore, if the condition \eqref{4:4}
holds then the converse is also true.
\et
\bn \label{rem4.1}
Observe that if, for instance, $|a|=\infty$
(the case $|b|=\infty$ is analogous) the
relations regarding to $a$ into \eqref{4:4}
need to be replaced by
$$
\lim_{t\to -\infty} x^k(t)x^l_{-1}(t)
\sigma(t)\rho(t)=0, \qquad  k,l\ge 0.
$$
\en
\bdm
Let $(p_n)$ be a
$q$-classical OPS with respect to
$\rho(s)$ on $x(s)$.
From  Lemma \ref{lem4.2} the function
$Q_k(x(s)):={\mathscr R}_1(\rho(s),x(s))
\big(p_{k-1}(x_1(s))\big)$ is a polynomial
on $x(s)$ of degree $k$.
Therefore for any $n>k\ge 1$
$$ \ba{rl}
0=&\displaystyle\sum_{s=a}^{b-1} p_n(s)Q_k(s)
\rho(s)\nabla x_1(s)\\[5mm] =&\displaystyle
\sum_{s=a}^{b-1} p_n(s)\Big({\mathscr R}_1
(\rho(s),x(s))\big(p_{k-1}(x_1(s))\big)
\Big)\,\rho(s)\nabla x_1(s) \\[5mm] =&
\displaystyle \sum_{s=a}^{b-1}p_n(s) \nabla
\Big(\rho_1(s) p_{k-1}(x_1(s))\Big).
\ea
$$
Taking into account the Leibniz rule, i.e.
$$
\nabla\big(f(s)g(s)\big)=\big(\nabla
f(s)\big)g(s)+f(s-1)\big(\nabla g(s)\big),
$$
the following formula holds
\begin{equation} \label{4:5}
\sum_{s=a}^{b-1}
g(s)(\nabla f(s))=f(s)
g(s+1)\Big|_{s=a-1}^{s=b-1}-\sum_{s=a}^{b-1}f(s)(\Delta g(s))
\end{equation}
Thus we get
$$
\ba{rl}
0= & \displaystyle p_n(s+1) \rho_1(s)
p_{k-1}(x_1(s))\Big|_{s=a-1}^{s=b-1}-
\sum_{s=a}^{b-1} \big(\Delta
p_n(s)\big)\rho_1(s)p_{k-1}(x_1(s))
\\[4mm]= & \displaystyle p_n(s) \rho(s)
\sigma(s)p_{k-1}(x_{-1}(s))\Big|_{s=a}^{s=b}
-\sum_{s=a}^{b-1}\big(\Delta p_n(s)\big)
\rho_1(s)p_{k-1}(x_1(s)) \quad (\Delta
x(s)=\nabla x_2(s))\\[4mm] = & \displaystyle
p_n(s)\sigma(s) \rho(s) p_{k-1}(x_{-1}(s))
\Big|_{s=a}^{s=b}-\sum_{s=a}^{b-1}\big(
\Delta^{(1)}p_n(s)\big) p_{k-1}(x_1(s))
\,\rho_1(s)\nabla x_2(s).
\ea
$$
Hence $(\Delta^{(1)}p_{n+1})$ is an OPS
with respect to $\rho_1(s)$ on the
lattice $x_1(s)$.
\\
Now, we prove that $\rho_1(s)$ satisfies
a $q$-Pearson equation as \eqref{4:2} on
the lattice $x_1(s)$.
Indeed since $(p_n)$ is a classical OPS
we get
$$\ba{rl}
\displaystyle \frac{\Delta (\sigma(s)
\rho_1(s))}{\rho_1(s)}= & \displaystyle
\frac{\sigma(s+1)\rho_1(s+1)}{\rho_1(s)}
-\sigma(s)\\ = & \displaystyle \sigma(s+
1)+\tau(x(s+1))\nabla x_1(s+1)-\sigma(s)
\\[4mm] =& \displaystyle \widehat{\sigma}
(s+1)-\widehat{\sigma}(s)+\frac{1}{2}
\tau(s+1)\nabla x_1(s+1)+ \frac{1}{2}
\tau(s)\nabla x_1(s),
\ea
$$
where $\widehat{\sigma}(s)=\sigma(s)+\frac
12 \tau(x(s))\nabla x_1(s)$ is a polynomial
on $x(s)$ of degree at most 2, and after a
straightforward calculation we deduce that
this expression is equal to $\widehat{\tau}_1
(x_1(s))\nabla x_2(s)$ with $\deg(\widehat
{\tau_1})\le 1$.
\\
Moreover, using the last expressions we get
$$
\sigma(s)+\frac{1}{2} \widehat{\tau}_1(x_1(s)
)\nabla x_2(s)=\frac{1}{2} \left(\widehat{
\sigma}(s+1)+\widehat{\sigma}(s)+\frac{1}{2}
\tau(s+1)\nabla x_1(s+1) -\frac{1}{2}\tau(s)
\nabla x_1(s) \right),
$$
is a polynomial on $x_1(s)$ of degree, at
most, 2 and hence the result holds.
\\
For the converse, we know that there exist
two polynomials on $x_1(s)$, $\hat \sigma(s)$
and $\tau_1$, of degree at most 2 and 1,
respectively, such that
$$
\Delta [\sigma(s)\rho_1(s)]=\tau_1(x_1(s))
\rho_1(s)\nabla x_2(s),
$$
here $\sigma(s):=\hat \sigma(s)-\frac{1}{2}
\tau_1(x_1(s))\nabla x_2(s)$.
\\
So, we only need to check that
$$
\tau(s):=\frac{\nabla x(s)}{\nabla x_1(s)}
\left(\tau_1(x_1(s-1))- \frac{\nabla \sigma(s)}
{\nabla x(s)}\right),
$$
and $\sigma(s)+\frac{1}{2} \tau(s)\nabla x_1(s)$
are polynomials on $x(s)$ of degree at
most 1 and 2, respectively, which is a direct
calculation and hence the result follows.
\edm
\bn \label{rem4.2}
Notice that in the case $a=-\infty$ (the case
$b=\infty$ is analogous) one should write
the formula \eqref{4:5} as
$$
\sum_{s=-\infty}^{b-1}
g(s)\nabla\big(f(s)\big)=f(b-1)g(b)-
\left(\lim_{t\to -\infty}f(t)g(t+1)
\right)-\sum_{s=-\infty}^{b-1} f(s)\Delta\big(g(s)
\big).
$$
\en
\bt \label{theo4.2}
Let $(p_n)$ be a $q$-classical OPS with
respect to $\rho(s)$ on $x(s)$ such that
\begin{equation} \label{4:6}
x^k(a)x^l_{-1}(a)\rho(a)=x^k(b)x^l_{-1}(b)
\rho(b)=0, \quad \forall\, k,l=0,1,\dots
\end{equation}
If $p_{-1}=0$, then  the sequence
$({\mathscr R}_n(\rho_{-1}(s),x_{-1}(s))(1))$
is a $q$-classical OPS with respect to the
weight function $\rho_{-1}(s)=\rho(s-1)/
\sigma(s)$ on $x_{-1}(s)$.
Furthermore, if the condition \eqref{4:4}
holds then the converse is also true.
\et
The proof follows the same steps as in
Theorem \ref{theo4.1} taking into account
the properties of $\Delta^{(m)}$ and the
basic relations of $\rho_m(s)$ and $x_m(s)$
for every integer $m$.
\bn \label{rem4.3}
\begin{enumerate} 
\item[(i)] The relation between the statements
of Theorems \ref{theo4.1} and \ref{theo4.2}
follows from
$$
{\mathscr R}_1(\rho_{-1}(s),x_{-1}(s))
\circ {\mathscr R}_n(\rho(s),x(s))=
{\mathscr R}_{n+1}(\rho_{-1}(s),x_{-1}(s)),
\quad n=0,1,2,\dots
$$
\item[(ii)] If the property of orthogonality
$$
\sum_{s=a}^{b-1} p_n(x(s))p_m(x(s))\rho(s)
\nabla x_1(s)=0,
$$
and the boundary conditions \eqref{4:4}
holds, then for every integer $k$ we get
$$
\sum_{s=a-k}^{b-k-1} \Big(\Delta^{(k)}
p_n(x_k(s))\Big)\Big(\Delta^{(k)}
p_m(x_k(s))\Big)\rho_k(s)\nabla x_{k+1}(s)=0,
$$
where if $k\ge 0$
$$
\Delta^{(-k)}:={\mathscr R}_k(\rho_{-k}(s),
x_{-k}(s)).
$$
\end{enumerate}
\en
Now we can state the first main result of
this paper.
\bt \label{theo4.4}
Let $(p_n)$ be an OPS with respect
to $\rho(s)$ on the lattice $x(s)$ defined in
\eqref{4:1} and let $\sigma(s)$ be such that
\eqref{4:6} holds.
Then the following statements are equivalent.
\begin{enumerate}
\item $(p_n)$ is $q$-classical.
\item The sequence $(\Delta^{(1)}p_n)$
is an OPS with respect to the weight function
$\rho_1(s)=\sigma(s+1)\rho(s+1)$
where $\rho$ satisfies \eqref{4:2}.
\item For every integer $k$, the sequence
$({\mathscr R}_n(\rho_k(s),x_k(s))(1))$ is
an OPS with respect to the weight function
$\rho_k(s)$ where $\rho_0(s)=\rho(s)$,
$\rho_k(s)=\rho_{k-1}(s+1)\sigma(s+1)$,
and $\rho$ satisfies \eqref{4:2}.
\item (Second order linear difference equation):
\label{4:12.1} $(p_n)$ satisfies the following
second order linear difference equation of hypergeometric type
\bq \label{4:18}
\sigma(s)\frac{\Delta}{\nabla x_1(s)}\frac
{\nabla p_n(s)}{\nabla x(s)} +\tau(s)\frac
{\Delta p_n(s)}{\Delta x(s)}+\lambda_n p_n(s)=0,
\eq
where $\widehat{\sigma}(s)=\sigma(s)+\frac{1}{2}
\tau(s)\nabla x_1(s)$ and $\tau(s)$ are
polynomials on $x(s)$ of degree at most 2
and 1, respectively, and $\lambda_n$ is
a constant.
\item $(p_n)$ can be expressed
in terms of the Rodrigues Operator as follows
\begin{equation} \label{4:13}
p_n(s)=B_n {\mathscr R}_n(\rho(s),x(s))(1)=
\frac{B_n}{\rho(s)}\frac{\nabla}{\nabla x_1(s)}
\frac{\nabla}{\nabla x_2(s)}\dots
\frac{\nabla}{\nabla x_n(s)}\,\left(
\rho_n(s)\right),
\end{equation}
where $B_n$ is a non zero constant.
\item \label{re7} ({\sf Second structure relation}):
There exist three sequences of complex numbers,
$(e_n)$, $(f_n)$, and $(g_n)$, such
that the fo\-llo\-wing relation holds
for every $n\ge 0$, with the convention $p_{-1}=0$,
$$
{\mathscr M} p_{n}(x(s))=e_n\frac{\Delta p_{n+1}(s)}
{\Delta x(s)}+ f_n \frac{\Delta p_{n}(s)}{\Delta x(s)}
+g_n\frac{\Delta p_{n-1}(s)}{\Delta x(s)},
$$
where $\mathscr M$ is the forward arithmetic mean
operator:
$$
{\mathscr M} f(s):=\frac{f(s+1)+f(s)}2,
$$
$e_n\ne 0$, $g_n\ne \gamma_n$ for all
$n\ge 0$, and $\gamma_n$ is the corresponding
coefficient of the three-term recurrence
relation \cite{nisuuv}
$$
x(s)p_n(s)=\alpha_n p_{n+1}(s)+
\beta_n p_n(s)+\gamma_n p_{n-1}(s),
\quad n\ge 0.
$$
\end{enumerate}
\et

\bn \label{rem4.4}
If we consider a $q$-linear lattice, i.e.
either $c_1=0$ or $c_2=0$, this result is
``equivalent'' to Theorem 1.3 stated
in \cite{alv3}, because $\tau(s)\nabla
x_1(s)$ is a polynomial on $x(s)$ of degree
2 and $x_k(s)=q^{\alpha_k}x(s)+\delta_k$,
where $\alpha_k,\,\delta_k$ are independent
of $s$ for every integer $k$.
\en

Assuming the theorem proved let us state and
proof the second main result.
\bt \label{theo4.3}
Under the hypothesis of Theorem \ref{theo4.4}
the following statements are equivalent.
\begin{enumerate}
\item[(i)] $(p_n)$ is $q$-classical.
\item[(ii)] $(\Delta^{(1)}p_{n+1})$ is a OPS.
\end{enumerate}
\et
\bdm
Of course (i)$\Rightarrow$(ii) is a consequence of
Theorem \ref{theo4.1}.
\\
(ii)$\Rightarrow$(i): This proof is analogous to the
proof given by W. Hahn in \cite{Hah}. We know that
the SPO $(p_n)$ satisfies the TTRR
$$
p_n(x(s))=(x(s)-a_n)p_{n-1}(x(s))-b_np_{n-2}(x(s)),
\quad n\ge 1,
$$
with initial conditions $p_{-1}=0$ and $p_0=1$.
Thus applying the (forward) arithmetic mean process
we get
\bq \label{4:7}
\Delta^{(1)} p_n(x(s))={\mathscr M} p_{n-1}(x(s))-b_n
\Delta^{(1)} p_{n-2}(x(s))+(\alpha_q(2)x_1(s)-\tilde
a_n)\Delta^{(1)} p_{n-1}(x(s)),
\eq
with $\tilde a_n=a_n-c_3(q)(1-\alpha_q(2))$.

Moreover, by hypothesis, the
sequence of monic polynomials $([n]^{-1}_q \Delta^{(1)}
p_n)$ also satisfies a TTRR, i.e.
\bq \label{4:8}
\frac 1{[n]_q} \Delta^{(1)} p_n(x(s))=(x_1(s)-a_n')
\frac 1{[n-1]_q}\Delta^{(1)} p_{n-1}(x(s))-\frac{b_n'}
{[n-2]_q}\Delta^{(1)} p_{n-2}(x(s)),
\quad n\ge 1.
\eq
Combining \eqref{4:7} and \eqref{4:8} to
eliminate $\Delta^{(1)}p_{n-2}$ (resp.
$\Delta^{(1)}p_{n}$)  we get
\bq \label{4:9} \begin{array}{rl}
\left(\frac{b_n'}{[n-2]_q}-\frac{b_n}{[n]_q}\right)
\Delta^{(1)}p_n(x(s))=& \left(\Big(\frac{\alpha_q(2)
b_n'}{[n-2]_{q}}-\frac{b_n}{[n-1]_q}\Big)x_1(s)-
\frac{b_n'\tilde a_n}{[n-2]_{q}}+\frac{b_na_n'}
{[n-1]_q}\right)\Delta^{(1)}p_{n-1}(x(s))\\[2mm] &
+\frac{b_n'}{[n-2]_{q}}{\mathscr M}p_{n-1}(x(s)),
\end{array}\eq
and
$$
\begin{array}{rl}
\frac{1}{[n]_q}{\mathscr M}p_{n-1}(x(s))=& \left(
\frac{[n]_q-[n-2]_q}{2[n-1]_q[n]_q}x_1(s)-\frac{a_n'}
{[n-1]_q}+\frac{\tilde a_n}{[n]_q}\right)\Delta^{(1)}
p_{n-1}(x(s))\\[2mm] & -\left(\frac{b_n'}{[n-2]_{q}}-
\frac{b_n}{[n]_q}\right)\Delta^{(1)}p_{n-2}(x(s)).
\end{array}
$$
Setting $n+1$ instead of $n$ in the above expression
the following expression fulfills:
\bq\label{4:10}
\begin{array}{rl}
\frac{1}{[n+1]_q}{\mathscr M}p_{n}(x(s))=& \left(
\frac{[n+1]_q-[n-1]_q}{2[n]_q[n+1]_q}x_1(s)-
\frac{a_{n+1}'}{[n]_q}+\frac{\tilde a_{n+1}}{[n+1]_q}
\right)\Delta^{(1)}p_{n}(x(s))\\[2mm] & +\left(
\frac{b_{n+1}}{[n+1]_q}-\frac{b_{n+1}'}
{[n-1]_{q}}\right)\Delta^{(1)}p_{n-1}(x(s)).
\end{array}\eq
In order to simplify further calculations in this
proof we denote by
$e_n:=\frac{b_n'}{[n-2]_{q}}-\frac{b_n}{[n]_q}$,
$k(x;n)=\frac{[n+1]_q-[n-1]_q}{2[n_q][n+1]_q}x-
\frac{a_{n+1}'}{[n]}+\frac{\tilde a_{n+1}}
{[n+1]_q}$ and by $l(x;n)=\Big(\frac{\alpha_q(2)
b_n'}{[n-2]_{q}}-\frac{b_n}{[n-1]_q}\Big)x-
\frac{b_n'\tilde a_n}{[n-2]_{q}}+\frac{b_na_n'}
{[n-1]_q}$.

Taking into account \eqref{4:9} and
\eqref{4:10} in order to eliminate $\Delta^{(1)}p_n$
we obtain
{\small \bq \label{4:11}
\frac {e_n}{[n+1]_q}{\mathscr M}p_n(x(s))=
\frac{b_n'}{[n-2]_{q}}k(x_1(s);n){\mathscr
M}p_{n-1}(x(s))+k(x_1(s);n)(l(x_1(s);n)-
e_{n+1}e_n)\Delta^{(1)}p_{n-1}(x(s)).
\eq}

By using \eqref{4:10} and \eqref{4:11} to
eliminate $p_{n}$ we get
\bq \label{4:11.1}
\begin{array}{c}
\frac{{\mathscr M}\big(l(x_{-1}(s);n)
\nabla p_{n-1}(x(s))\big)}{\nabla x_1(s)}
+\frac{b_n'}{[n-2]_q}\frac{{\mathscr
M}\big((\nabla x(s)){\mathscr M}p_{n-1}
(x(s-1))\big)}{\nabla x_1(s)} \\[3mm]
=\frac{[n+1]_q b_n'}{[n-2]_q}\frac{\nabla
\left(k(x_1(s);n){\mathscr M}p_{n-1}(x(s))
\right)}{\nabla x_1(s)}+\mbox{\scriptsize
$[n\!+\!1]_q$}\frac{\nabla\left((l(x_1(s);n)
k(x_1(s);n)-e_ne_{n+1})\Delta^{(1)}
p_{n-1}(x(s))\right)}{\nabla x_1(s)}.
\end{array}\eq

On the other hand, combining \eqref{4:9} and
\eqref{4:10} to eliminate $\Delta^{(1)}
p_{n-1}$ we get
\bq \label{4:12}\begin{array}{c}
\frac{b_n'e_{n+1}}{[n-2]_q}{\mathscr M}
p_{n-1}(x(s))=\frac{l(x_1(s);n)}{[n+1]_q}
{\mathscr M}p_n(x(s))-(k(x_1(s);n)
l(x_1(s);n)-e_ne_{n+1})\Delta^{(1)}p_n(x(s)).
\end{array}\eq
By using the expressions \eqref{4:10} and
\eqref{4:12} to eliminate $p_{n-1}$ we get
\bq \label{4:11.2}
\begin{array}{c}
\frac{1}{[n+1]_q}\frac{{\mathscr M}
\big((\nabla x(s)){\mathscr M}p_{n}
(x(s-1))\big)}{\nabla x_1(s)}-\frac{
{\mathscr M}\big(k(x_{-1}(s);n)\nabla
p_{n}(x(s))\big)}{\nabla x_1(s)}\\[3mm]
=\frac{[n-2]_q}{b_n'[n+1]_q}\frac{\nabla
\left(l(x_1(s);n){\mathscr M}p_{n}(x(s))
\right)}{\nabla x_1(s)}-\frac{[n-2]_q}
{b_n'}\frac{\nabla\left((l(x_1(s);n)
k(x_1(s);n)-e_ne_{n+1})\Delta^{(1)}
p_{n}(x(s))\right)}{\nabla x_1(s)}.
\end{array}\eq

The following result has a technical role and it will be used in to proof of this Theorem.

\bl \label{lem4.3} For any two polynomials
$\pi$ and $\xi$, the following relation holds:
$$
\begin{array}{rl}
\frac{{\mathscr M}\left(\pi(x_{-1}(s))
\nabla \xi (x(s))\right)}{\nabla x_1(s)}
=& \frac{\nabla(\pi(x_{1}(s)){\mathscr M}
\xi(x(s)))}{\nabla x_1(s)}+(\nabla^{(1)}
\pi(x_1(s))) \xi(x(s))\\[2mm] =&\frac{
{\mathscr M}(\pi(x_{-1}(s))\nabla x(s))}
{\nabla x_1(s)}{\mathscr M}\frac{\nabla
\xi(x(s))}{\nabla x(s)}+\nabla\big(\pi(x_1(s))
\Delta x(s)\big)\frac{\nabla \Delta^{(1)}
\xi(x(s))}{\nabla x_1(s)}.
\end{array}$$
Moreover, the functions
$\frac{{\mathscr M}(\pi(x_{-1}(s))\nabla x(s))}
{\nabla x_1(s)}$ and $\nabla\big(\pi(x_1(s))
\Delta x(s)\big)$ are polynomials on $x(s)$.
\el

The proof of this result is
straightforward and we leave it as an exercise
for the reader.

Therefore, using Lemma
\ref{lem4.3} and  doing some extra
calculations, we obtain that the
expressions \eqref{4:11.1} and
\eqref{4:11.2} become the
following two second order homogeneous
linear difference equations:
\begin{eqnarray}
\nonumber \phi_2(x(s),n)\nabla^{(1)}
\Delta^{(1)}p_n(x(s))+\phi_1(x(s),n)
{\mathscr M}\frac{\nabla p_n(x(s))}
{\nabla x(s)}+\phi_0(n) p_n(x(s))=0,
\hspace{11mm} \\ \nonumber \phi_2(x
(s),n)\nabla^{(1)}\Delta^{(1)}p_{n-1}
(x(s))+\phi_1(x(s),n){\mathscr M}
\frac{\nabla p_{n-1}(x(s))}{\nabla x(s)}+
\tilde \phi_0(n) p_{n-1}(x(s))=0,
\end{eqnarray}
where the indices indicate the degree
of the polynomial coefficients.

Moreover, replacing $n$ by $n+1$ in the
last expression we obtain the second
order linear difference equation:
\begin{equation}
\phi_2(x(s),n+1)\nabla^{(1)}\Delta^{(1)}
p_{n}(x(s))+ \phi_1(x(s),n+1){\mathscr M}
\frac{\nabla p_{n}(x(s))}{\nabla x(s)}+
\tilde \phi_0(n+1) p_{n}(x(s))=0.
\end{equation}
So, we have two difference equations
for $p_n$ which only differ in a constant
factor.

Therefore $p_n$ satisfies the following
two second order linear
difference equations:
$$
\varphi_2(x(s))\nabla^{(1)}\Delta^{(1)}
p_n(x(s))+\varphi_1(x(s)){\mathscr M}
\frac{\nabla p_{n}(x(s))}{\nabla x(s)}
+\varphi_0(x(s))p_{n}(x(s))=0,
$$
$$
\psi_2(x(s))\nabla^{(1)}\Delta^{(1)}p_n(x(s))+
\psi_1(x(s)){\mathscr M}\frac{\nabla p_{n}(x(s))}
{\nabla x(s)}+\psi_0(x(s))p_{n}(x(s))=0,
$$
where, again, the indices indicate the
degree of the polynomial coefficients.
Denoting by $\sigma(s):=\varphi_2(x(s))+
\varphi_1(x(s))\nabla x_1(s)$, $\zeta(s):=
\psi_2(x(s))+\psi_1(x(s))\nabla x_1(s)$, the
above expressions can be rewritten as
$$
\sigma(s)\nabla^{(1)}\Delta^{(1)}
p_n(x(s))+\varphi_1(x(s))\Delta^{(1)} p_n(x(s))
+\varphi_0(n)p_{n}(x(s))=0,
$$
$$
\zeta(s)\nabla^{(1)}\Delta^{(1)}
p_n(x(s))+\psi_1(x(s))\Delta^{(1)} p_n(x(s))
+\psi_0(n)p_{n}(x(s))=0.
$$
Notice that $\sigma$ and $\zeta$ are not
always polynomials on $x(s)$.

So,
$$\ba{l}
(\varphi_1(x(s))\zeta(s)-\sigma(s)\psi_1(x(s)))
\Delta^{(1)} p_{n}(x(s))\\[1mm] +(\varphi_0(n)
\zeta(s)-\sigma(s)\psi_0(n))p_n(x(s))=0.
\ea$$
If $n>3$, we get that $p_n$ and $\Delta^{(1)}
p_n$, has a common zero but $p_n$ has not
no multiple zeros.
Hence, after dividing by the common
factor, the polynomial coefficients
of the difference equations for  $p_{n+1}$,
$p_n$, and $p_{n-1} $ are the same.
Taking this into account follow the coefficients do not depend on
$n$, so $p_n(x)$ satisfies the second
order linear difference
equation
$$
\phi(x(s))\nabla^{(1)}\Delta^{(1)}
p_n(x(s))+\psi(x(s)){\mathscr M}
\frac{\nabla p_n(x(s))}{\nabla x(s)}+
\lambda_n p_n(x(s))=0,
$$
where $\deg \phi\le 2$, $\deg \psi\le
1$, and  $\lambda_n=-[n]_q
(\alpha_q(n-1)b_\psi+[n-1]a_\phi)$
depends on $n$ since $p_n$ is a
polynomial, being $a_\phi$ and
$b_\psi$ the leading coefficients of
$\phi$ and $\psi$, respectively.

Hence, according to the \ref{4:12.1}th condition of
Theorem \ref{theo4.4}, $(p_n)$ is $q$-classical.
\edm


Let us prove now Theorem \ref{theo4.4}.
From Theorem \ref{theo4.1}, Theorem \ref{theo4.2},
and its Corollary, as well as Remark \ref{rem4.3}
we know that $(1)\rightarrow (2)\rightarrow
(3)\rightarrow (1)$.
\bp \label{prop4.1}{\rm ((1)$\rightarrow$ (4))}
If $(p_n)$ is a $q$-classical
OPS with respect to $\rho(s)$ satisfying the
boundary relations \eqref{4:6} then $(p_n)$
satisfies the second order
linear difference equation of hypergeometric type
\eqref{4:18}.
\ep
\bdm
By Lemma \ref{lem4.2} the function
${\mathscr R}_1(\rho(s),x(s))
\big(\Delta^{(1)}p_n(x(s)\big)$ is a
polynomial on $x(s)$ of degree $n$.
On the other hand, by Theorem \ref{theo4.1}
$(\Delta^{(1)}p_n)$ is a $q$-classical
OPS with respect to $\rho_1(s)$
on $x_1(s)$. Thus, for $0\le k<n$,
$$\ba{rl}
& \displaystyle \sum_{s=a}^{b-1}Q_k(x(s))
{\mathscr R}_1(\rho(s),x(s))\big(\Delta^{(
1)}p_n(x(s))\big)\rho(s)\nabla x_1(s) \\[0.4cm]
=&\displaystyle \sum_{s=a}^{b-1} Q_k(x(s))
\nabla\big(\rho_1(s)\Delta^{(1)}p_n(x(s))
\big) \\[0.4cm] \stackrel{\eqref{4:5}}= &
\displaystyle Q_k(x(s))\sigma(s) \rho(s)
\Delta^{(1)} p_n(x(s-1))\Big|_{s=a}^{s=b}
\!\!- \!\! \sum_{s=a}^{b-1} \Delta^{(1)}
\big(p_n(x(s))\big)\Delta^{(1)}
\big(Q_k(x(s))\big)\, \rho_1(s) \nabla x_2(s).
\ea
$$
By hypothesis and Theorem \ref{theo4.1}
the last expression vanishes and, as a
consequence, there exists a non zero
constant $\lambda_n$, independent of
$s$, such that
$$
{\mathscr R}_1(\rho(s),x(s))
\big(\Delta^{(1)}p_n(x(s))\big)=
-\lambda_n p_n(x(s)),
$$
since $(p_n)$  is an OPS with respect
to $\rho(s)$ on $x(s)$.
\\
Finally, using that the above expression is
indeed ${\mathscr H}_q$ (see \eqref{aux1})
it is clear that the above expression
\eqref{4:18} and hence the result holds.
\edm
\\
(4)$\rightarrow$(5):
This is a well-known property.
See, for instance, \S 3.2.2 in page 64 of
\cite{nisuuv}.
\\
(5)$\rightarrow$(1): By setting $n=1$ in the
Rodrigues formula we obtain the Pearson type
equation \eqref{4:2} that the
weight function $\rho(s)$ satisfies.

\bp \label{prop4.2}$(\eqref{re7}\leftrightarrow(1))$.
Let $(p_n)$ be an OPS with respect
to $\rho(s)$ on the lattice $x(s)$ defined
in \eqref{4:1} and let $\sigma(s)$ be such
that \eqref{4:6} holds.
The sequence $(p_n)$ is $q$-classical
if and only if there exist sequences of complex
numbers, $(e_n)$, $(f_n)$, and $(g_n)$,
such that the following relation holds
$$
{\mathscr M}p_n(x(s))=e_n
\frac{\Delta p_{n+1}(x(s))}{\Delta x(s)}+
f_n \frac{\Delta p_{n}(x(s))}{\Delta x(s)}
+g_n\frac{\Delta p_{n-1}(x(s))}{\Delta x(s)},
$$
with the convention $p_{-1}=0$,
where $e_n\ne 0$, $g_n\ne \gamma_n$ for
all $n\ge 0$.
\ep
\bdm If $(p_n)$ is an OPS then it
satisfies a TTRR, i.e. there exist three
sequences of complex numbers, $(\alpha_n)$,
$(\beta_n)$, and $(\gamma_n)$, with
$\gamma_n\ne 0$, such that
\begin{equation} \label{4:15}
x(s)p_n(x(s))=\alpha_n p_{n+1}(x(s))+
\beta_n p_n(x(s))+\gamma_n p_{n-1}(x(s)),
\quad n\ge 0.
\end{equation}
If $(p_n)$ is $q$-classical, then
$(\Delta^{(1)}p_{n+1})$ is an OPS,
and therefore there exist sequences of
complex numbers, $(\alpha^{(1)}_n)$,
$(\beta^{(1)}_n)$, and $(\gamma^{(1)}_n)$,
with $\gamma^{(1)}_n\ne 0$,
such that
\begin{equation} \label{4:16}
x_1(s)\Delta^{(1)}p_n(x(s))=
\alpha^{(1)}_n \Delta^{(1)}
p_{n+1}(x(s))+\beta^{(1)}_n
\Delta^{(1)}p_n(x(s))+ \gamma^{(1)}_n
\Delta^{(1)}p_{n-1}(x(s)).
\end{equation}
But
\begin{equation} \label{4:17}
\Delta^{(1)}\big(x(s)p_n(x(s))\big)=
({\mathscr M} x(s))\Delta^{(1)}p_n(x(s))+
{\mathscr M} p_n(x(s)).
\end{equation}
So combining \eqref{4:15}, \eqref{4:16},
and \eqref{4:17} we get
$$
{\mathscr M} p_n(x(s))\!=\!\left(\!\alpha_n-
\frac{[2]_q}{2}\, \alpha^{(1)}_n\!\right)
 \!\!\Delta^{(1)} p_{n+1}(x(s))+ \widehat
 \beta_n \Delta^{(1)}p_n(x(s))+
\left(\!\gamma_n-\frac{[2]_q}2
\gamma^{(1)}_n\!\right)\!\!
\Delta^{(1)}p_{n-1}(x(s)),
$$
where
$$
\widehat\beta_n=(\beta_n-c_3)
-\frac{[2]_q}2(\beta^{(1)}_n-c_3).
$$
Moreover, the coefficient of
$\Delta^{(1)}p_{n-1}$ is different
from $\gamma_n$ because $\gamma^{(1)}_n\ne 0$.
\\
Conversely, if there exist sequences of
complex numbers, $(e_n)$, $(f_n)$, and
$(g_n)$, such that the following relation
holds
$$
{\mathscr M} p_n(x(s))=e_n\frac{\Delta
p_{n+1}(x(s))}{\Delta x(s)}+f_n \frac
{\Delta p_{n}(x(s))}{\Delta x(s)}
+g_n\frac{\Delta p_{n-1}(x(s))}
{\Delta x(s)},\qquad n\ge 0,
$$
assuming $p_{-1}=0$, then from
\eqref{4:17} we get
$$
\frac{[2]_q}2x_1(s)\Delta^{(1)}
p_n(x(s))=(\alpha_n-e_n)\Delta^{(1)}
p_{n+1}(x(s))+ \widehat \beta^{(1)}_n
\Delta^{(1)}p_n(x(s))+(\gamma_n-g_n)
\Delta^{(1)}p_{n-1}(x(s)),
$$
where
$$
\widehat \beta^{(1)}_n=
\beta_n-c_3-f_n+\frac{[2]_q}2c_3.
$$
By hypothesis $g_n\ne \gamma_n$, hence by
Favard's Theorem $(\Delta^{(1)}p_{n+1})$
is an OPS, and by Hahn's Theorem
 the result holds.
\edm
\section{The examples} \label{sec5}
\subsection{The Askey-Wilson Polynomials}
\label{sec5.1}
The Askey-Wilson polynomials, which
were introduced by R. Askey and J.
Wilson, are located in the top of
the $q$-Askey  tableau \cite{kost}.
These polynomials can be written as
a basic hypergeometric series
$$
p_n(x(s);a,b,c,d|q)=\frac{(ab;q)_n
(ac;q)_n(ad;q)_n}{a^n}\,_4\varphi_3
\left(\ba{c|}q^{-n},\, abcdq^{n-1},\,
aq^s, \, aq^{-s} \\ ab,\, ac,\, ad
\ea\ q;q\right),
$$
where $a, b, c, d, ab, ac, ad, bc,
bd, cd\not \in \{q^{m}:m\in {\mathbb Z}\}$.
\bn \label{rem51}
Notice that although this family does not satisfy
a property of orthogonality \eqref{2:10} it is
orthogonal with respect the linear functional
${\bf u}^{AW}$ which has the following integral
representation \cite[(3.1.2), p. 63]{kost}
$$
\pe {{\bf u}^{AW}}p=
\int_{-1}^1 p(x)\widetilde \omega(x)dx=
\int_{-1}^1 \frac{p(x)}{\sqrt{1-x^2}}
\left|\frac{(e^{2i\theta};q)_\infty}
{(ae^{i\theta},be^{i\theta},ce^{i\theta},
de^{i\theta};q)_\infty}\right|^2dx, \quad x=\cos\theta.
$$
Furthermore, such a functional fulfills the
distributional equation
$$
\Delta^{(1)}(\sigma^{AW}{\bf u}^{AW})=
p_1(x(s);a,b,c,d|q){\bf u}^{AW},
$$
and it is straightforward to check that
$(\sigma^{AW}(s),\omega(s),x(s))$ is
$q$-classical, where $x(s)\!=\!\frac{1}{2}\!
\left(q^s\!+\!q^{-s}\right)$, i.e. $c_1=c_2=\frac{1}
{2}$ and $c_3=0$, $\sigma^{AW}(s)=
-\kappa_q^2q^{-2s+\frac{1}{2}}(q^s-a)
(q^s-b)(q^s-c)(q^s-d)$,
and hence we can apply Theorem \ref{theo4.4}
to Askey-Wilson polynomials.
\en
These polynomials are the polynomial
eigenfunctions of the
linear difference operator \cite[(3.1.6)]{kost}
$$
{\mathscr H}^{AW}_q=\frac 1
{\nabla x_1(s)}\left(\sigma^{AW}
(-s)\frac{\Delta}{\Delta x(s)}-
\sigma^{AW}(s)\frac{\nabla}
{\nabla x(s)}\right).
$$
Notice that if $\tau^{AW}(x(s))=
p_1(x(s);a,b,c,d|q)$ a straightforward
calculation yields to
\begin{equation} \label{5:1}
\sigma^{AW}(-s)=\sigma^{AW}(s)+
\tau^{AW}(x(s))\nabla x_1(s).
\end{equation}

Here, as before, $\kappa_q=
q^\frac{1}{2}-q^{-\frac{1}{2}}$.
With these parameters, let us
consider the following function
$$
\rho^{AW}(s)=q^{-2s^2}(a,b,c,d;q)_s(a,b,c,d;q)_{-s},
$$
where $(a;q)_0=1$, $(a;q)_k=(1-a)(1-aq)
\cdots (1-aq^{k-1})$, $k\ge 1$,
and
$$
(a;q)_{-k}=\frac1{(1-aq^{-1})(1-aq^{-2})
\cdots (1-aq^{-k})}, \quad k\ge 1,
$$
\bl \label{lem5.1}
The function $\rho^{AW}(s)$ satisfies
the Pearson-type equation
\begin{equation} \label{5:2}
\sigma^{AW}(s+1)\rho(s+1)=\sigma^{AW}(-s)\rho(s).
\end{equation}
\el
\bdm
$$\ba{rl}
\displaystyle \frac{\rho^{AW}(s+1)}
{\rho^{AW}(s)}=& \displaystyle\frac
{q^{-2(s+1)^2}}{q^{-2s^2}}
\frac{(a,b,c,d;q)_{s+1}(a,b,c,d;q)_{-s-1}}
{(a,b,c,d;q)_{s}(a,b,c,d;q)_{-s}}\\[4mm]
=& \displaystyle\frac{q^{-4s-2}(1-aq^{s})
(1-bq^{s})(1-cq^{s})(1-dq^{s})}
{(1-aq^{-s-1})(1-bq^{-s-1})(1-cq^{-s-1})(1-
dq^{-s-1})}\\[4mm] = & \displaystyle
\frac{q^{2s}(q^{-s}-a) (q^{-s}-b)(q^{-s}-c)
(q^{-s}-d)}{q^{-2s-2}(q^{s+1}-a)(q^{s+1}
-b)(q^{s+1}-c)(q^{s+1}-d)}= \displaystyle
\frac{\sigma^{AW} (-s)}{\sigma^{AW}(s+1).}
\ea$$
\edm

Notice that from \eqref{5:1}, the
equation \eqref{5:2} becomes the
Pearson-type equation \eqref{4:2}.

Taking into account that $\rho_n(s)=
\rho_{n-1}(s+1) \sigma(s+1)$, $n\ge 1$,
with $\rho_0(s)=\rho(s)$, a
straightforward calculation gives
$$
\rho_n^{AW}(s)=\kappa_q^{2n}
q^{-2s^2-2sn-n^2+\frac{3}{2}n}
(a,b,c,d;q)_{s+n}(a,b,c,d;q)_{-s}.
$$
Then the Askey-Wilson polynomials
can be written for every
nonnegative integer $n$ as
$$
p_n(x(s);a,b,c,d|q)=\frac{B_n
\kappa_q^{2n}q^{2s^2}}{(a,b,c,d;q)_s
(a,b,c,d;q)_{-s}} \nabla^{(n)}
\frac{(a,b,c,d;q)_{s+n}(a,b,c,d;q)_{-s}}
{q^{2s^2+2sn+n^2-\frac{3}{2} n}},
$$
where $B_n=2^{-n}\kappa_q^{-n}
q^{\frac{n(3n-5)}4}$.
\\
Furthermore, since this family
fulfills the following difference
relation \cite[(3.1.8)]{kost}
$$
\Delta^{(1)}p_n(x(s);a,b,c,d|q)=
2[n]_q(1-abcdq^{n-1}) p_{n-1}
(x_1(s);aq^\frac{1}{2}, bq^\frac{1}{2},
cq^\frac{1}{2}, dq^\frac{1}{2}|q),
$$
and the coefficients of the second
structure relation are
$$\ba{rl}
e_n= & \displaystyle \frac{2[n]_q
(1-abcdq^{n-1})^2-[2]_q[n+
1]_q(1-abcdq^{n})^2}{4[n]_q(1-
abcdq^{n-1})(1-abcdq^{2n-1})
(1-abcdq^{2n})},\\[4mm] f_n=
& \displaystyle \frac{1-q}4(a-
a^{-1}q^{-1})-\left.\frac{1}{2}
\right(A_n(a,b,c,d|q)+C_n(a,b,c,d|q)
\\[4mm] & \displaystyle \left.-
\frac{[2]_q}2(A_{n-1}(aq^\frac{1}{2},
bq^\frac{1}{2},cq^\frac{1}{2},
dq^\frac{1}{2}|q)+C_{n-1}(a
q^\frac{1}{2},bq^\frac{1}{2},cq^\frac{1}{2}
,dq^\frac{1}{2}|q))\right),\\[4mm]
g_n= & \displaystyle \frac{(1-abcd
q^{[n-2]_{q}})(1-abcdq^{2n-2})(1-
abcdq^{2n-1})}{4[n]_q(1-abcdq^{n-1})^2}
\left(\frac{2[n]_q(A_{n-1}C_n)(a,b,c,d|q)}
{(1-abcdq^{[n-2]_{q}})^2}\right.\\[4mm]
& \displaystyle \left.-\frac{[2]_q[n-1]_q
(A_{[n-2]_{q}}C_{n-1})(aq^\frac{1}{2},
bq^\frac{1}{2}, cq^\frac{1}{2},d
q^\frac{1}{2}|q)}{(1-abcdq^{n-1})^2}
\right).
\ea
$$
\bn \label{rem5.1}
Notice that in \cite{koo4} T.H. Koornwinder
has obtained a different explicit structure
relation for Askey-Wilson polynomials by
using the difference operator theory.
\en
\subsection{The $q$-Racah polynomials} \label{sec5.2}
We consider the $q$-Racah polynomials
$u_n^{\alpha,\beta}(x(s),a,b)_q$
on the lattice
$$
x(s)=[s]_q[s+1]_q=q^\frac{1}{2}
\kappa_q^{-2} q^s+q^{-\frac{1}
{2}}\kappa_q^{-2}q^{-s}-[2]_q\kappa_q^{-2}.
$$
They were introduced in \cite{mal,nisuuv}
and deeply studied in \cite{alsmco}.
These polynomials can be written in terms of a basic
hypergeometric series \cite{alsmco} as follows
$$
u_n^{\alpha,\beta}(x(s),a,b)_q=
D_n \,_{4}\varphi_3 \left(\ba{c} q^{-n},
q^{ \alpha+\beta+n+1}, q^{a-s}, q^{a+s+1}
\\ q^{a-b+1},q^{\beta+1},q^{ a+b+\alpha+1}
\ea \,\bigg|\, q \,,\, q \right),
$$
where
$$
D_n=\frac{q^{-\frac{n}{2}(2a+\alpha+\beta+n+1)}
(q^{a-b+1};q)_n (q^{\beta+1};q)_n
(q^{a+b+\alpha+1};q)_n}{\kappa_q^{2n} (q;q)_n}.
$$

Observe that from the above formulas
the polynomials $u_n^{\alpha,\beta}
(x(s),a,b)_q$ are multiples of the
standard $q$-Racah polynomials
${\mathscr R}_n(\mu(q^{b+s});q^\alpha,
q^\beta,q^{a-b},q^{-a-b}|q)$.

These polynomials are the polynomial eigenfunctions
of the second order homogeneous linear difference
hypergeometric operator
$$
{\mathscr H}^{qR}_q=\frac 1{\nabla x_1(s)}
\left(\sigma^{qR}(-s-1)\frac{\Delta}{\Delta
x(s)}-\sigma^{qR}(s)\frac{\nabla}
{\nabla x(s)}\right),
$$
where
$$
\sigma^{qR}(s)=-\frac{q^{-2s}}{\kappa_q^4
q^{\frac{\alpha+\beta}2}}(q^s-q^a)(q^s-
q^{-b})(q^s-q^{\beta-a})(q^s-q^{b+
\alpha})=[s-a]_q[s+b]_q[s+a-\beta]_q
[b+\alpha-s]_q.
$$
Furthermore, taking into account the
$q$-difference relation
$$
\Delta^{(1)}u_n^{\alpha,\beta}
(x(s),a,b)_q=[\alpha+\beta+n
+1]_q u_{n-1}^{\alpha+1,\beta+1}
(x_1(s),a+\mbox{$\frac 12$},
b-\mbox{$\frac 12$})_q.
$$
The coefficients for the second
structure relation are
$$\ba{rl}
e_n= & \displaystyle \frac{2[n+\alpha+
\beta+1]_q-[2]_q[n+\alpha+\beta+2]_q}
{2[2n+\alpha+\beta+1]_q[2n+\alpha+\beta
+2]_q}[n+1]_q, \\[4mm] f_n= &
\displaystyle \beta_n(a,b,
\alpha,\beta)-\frac{[2]_q}2\beta_n
(a+\mbox{$\frac 12$},b-\mbox{$\frac 12$},
\alpha+1,\beta+1),\\[4mm] g_n= &
\displaystyle \frac{[a+b+\alpha+n]_q
[a+b-\beta-n]_q[\alpha+n]_q[\beta+n]_q
[b-a+\alpha+\beta+n]_q[b-a-n]_q}{2[\alpha+
\beta+2n]_q [\alpha+\beta+2n+1]_q[n+
\alpha+\beta+1]_q}\\[4mm] & \displaystyle
\times(2[n+\alpha+\beta+1]_q-[2]_q[n+\alpha+
\beta]_q).
\ea$$

\subsection{The $q$-Meixner polynomials}
\label{sec5.3}
This family of $q$-polynomials has
the following representation as a
basic hypergeometric series
$$
M_n(x;b,c;q)=\frac{(bq;q)_n(-c)^n}
{q^{n^2}}\,_2\varphi_1 \left.
\left(\ba{c}q^{-n},x\\ bq \ea
\right|q;-\frac{q^{n+1}}c\right),
\quad x(s)=q^{-s}\equiv x.
$$
They are the polynomial eigenfunctions
of the linear
difference operator of hypergeometric type
\cite[(3.13.5)]{kost}
$$
{\mathscr H}^{qM}_q=\frac 1{\nabla x_1(s)}
\left((x-1)(x+bc)\frac{\Delta}{\Delta x(s)}
-q^{-1}c(x-bq)\frac{\nabla}
{\nabla x(s)}\right),
$$
In this case
$$
\ba{rl}
e_n=& \displaystyle 1-\frac{[2]_q[n+1]_q}
{2[n]_q}q^{-\frac{1}{2}},\\[4mm] f_n=&
\displaystyle \frac{1-q}2+\frac{q^{-2n-1}}
2(c-cbq^{n+1})(2-q-q^2)+\frac{q^{-2n}}2
(c+q^n)(2-q-q^2+q^n-q^{n-1}),\\[4mm]
g_n=& \displaystyle \frac{cq^{-4n+1}(1-bq^n)
(c+q^n)}{2[n]_q}(2[n]_q(1-q^n)-
(q+q^2)[n-1]_q(1-q^{n-1})).
\ea
$$
\subsection{The Al-Salam \& Carlitz
polynomials I and II} \label{sec5.4}
The Al-Salam \& Carlitz polynomials I
(and II) appear in certain models of
$q$-harmonic oscillator, see e.g.
\cite{assu1,atsu2,atsu1,nag}.
They are polynomials orthogonal on the $q$-linear
lattice $x(s)=q^s\equiv x$, defined
\cite{kost} by
$$
U_n^{(a)}(x;q)=(-a)^n q^{n\choose 2}\, _2
\varphi_1\left(\begin{array}{c|c} q^{-n},
x^{-1} \\[-0.3cm] & q;\displaystyle
\frac{qx}{a} \\[-0.2cm]0 \end{array} \right).
$$
These polynomials are the polynomial
eigenfunctions of the linear difference operator of hypergeometric
type \cite[(3.24.5)]{kost}
$$
{\mathscr H}^{ACI}_q=\frac 1{\nabla
x_1(s)}\left(a \frac{\Delta}{\Delta
x(s)}-(x-1)(x-a)\frac{\nabla}
{\nabla x(s)}\right),
$$
In this case \cite[(3.24.7)]{kost}
$\Delta^{(1)}U_n^{(a)}(x;q)=
q^{\frac{n-1}2}[n]_qU_{n-1}^{(a)}(x;q)$,
 so the second structure relation coefficients are
$$\ba{rlrl}
e_n=& \displaystyle 1-\frac{(1+q)[n+1]_q}{2q[n]_q},&
f_n=& \displaystyle (1+a)q^n\left(1-\frac{[2]_q}
2q^{-\mbox{\scriptsize $\frac{3}{2}$} }\right),\\[4mm]
\multicolumn{4}{c}{g_n= \displaystyle \frac{aq^{n-\frac
52}}{2[n]_q}\left(2[n]q^{\mbox{$\frac{3}{2}$}}(q^n-1)-
[2]_q[n-1]_q(q^{n-1}-1)\right).}
\ea
$$
Taking into account that \cite[p. 115]{kost}
$$
V_n^{(a)}(x;q)=U_n^{(a)}(x;q^{-1}),
$$
all the information related to the Al-Salam \& Carlitz
polynomials II can be deduced from the information for
the Al-Salam \& Carlitz polynomials I in a simple way
by changing $q$ to $q^{-1}$.
\begin{acknowledgements}
The work of the authors has been supported by Direcci\'{o}n
General de Investigaci\'{o}n, Ministerio de Educaci\'{o}n y
Ciencia of Spain under grant MTM 2006-13000-C03-02.
Finally, we thank the anonymous referees for the careful
revision of the manuscript.
Their comments  and suggestions have  substantially improved the
presentation of the manuscript.
\end{acknowledgements}

\end{article}

\begin{thebibliography}{99}

\bibitem{alch}  
W. A. Al-Salam and T. S. Chihara.
\newblock Another characterization of the classical
orthogonal polynomials.
\newblock {\em SIAM J. Math. Anal.} {\bf 3} , 65--70 (1992)

\bibitem{alal} 
M. Alfaro and R. \'{A}lvarez-Nodarse.
\newblock A characterization of the classical orthogonal
discrete and $q$-polynomials.
\newblock{\em J. Comput. Appl. Math.} {\bf 201}, 48--54 (2007)

\bibitem{alv2} 
R. \'{A}lvarez-Nodarse.
\newblock {\em Polinomios hipergeom\'{e}tricos y
$q$-polinomios}, volume 26 of {\em Monograf{\'i}as del
Seminario Matem\'{a}tico ``Garcia de Galdeano''}.
\newblock Prensas Universitarias de Zaragoza, Zaragoza,
Spain. In Spanish (2003)

\bibitem{alv3} 
R. \'{A}lvarez-Nodarse.
\newblock On characterizations of classical polynomials.
\newblock {\em J. Comput. Appl. Math.} {\bf 196}, 320--337
(2006)

\bibitem{alar1} 
R. \'{A}lvarez-Nodarse and J. Arves\'{u}.
\newblock On the $q$-polynomials on the exponential lattice
$x(s)=c_1q^s+c_3$.
\newblock {\em Integral Trans. and Special Funct.}  {\bf 8},
299--324 (1999)

\bibitem{alsmco} 
R. \'{A}lvarez-Nodarse, Yu. F. Smirnov, and R. S.
Costas-Santos.
\newblock A $q$-analog of the Racah polynomials and $q$-algebra
$su_q(2)$ in quantum optics.
\newblock {\em Journal of Russian Laser Research}
{\bf 27}(1), 1--32 (2006)

\bibitem{and} 
G. E. Andrews.
\newblock {\em $q$-Series: Their Development and Application
in Analysis, Number Theory, Combinatorics, Physics and
Computer Algebra}, volume 66 of {\em C.B.M.S. Regional
Conference Series in Math}.
\newblock Amer. Math. Soc., Providence, Rhode Island (1996)

\bibitem{assu1} 
R. Askey and S. K. Suslov.
\newblock The $q$-harmonic oscillator and the Al-Salam and
Carlitz polynomials.
\newblock {\em Lett. Math. Phys.}  {\bf 29}(2), 123--132
(1993)

\bibitem{atsu2} 
N. M. Atakishiev and S. K. Suslov.
\newblock Difference analogs of the harmonic oscillator.
\newblock {\em Theoret. and Math. Phys.} {\bf 85}(1) ,
1055--1062 (1991)

\bibitem{atsu1} 
N. M. Atakishiev and S. K. Suslov.
\newblock A realization of the $q$-harmonic oscillator.
\newblock {\em Theoret. and Math. Phys.}{\bf 87}(1),
442--444 (1991)

\bibitem{atrasu} 
N. M. Atakishiyev, M. Rahman, and S. K.Suslov.
\newblock On classical orthogonal polynomials.
\newblock {\em Constr. Approx.} {\bf 11}, 181--226 (1995)

\bibitem{ch} 
T. S. Chihara.
\newblock {\em An Introduction to Orthogonal Polynomials},
\newblock  Gordon and Breach Science Publishers, New York
(1978)

\bibitem{fin} 
N. J. Fine.
\newblock {\em Hypergeometric Series and Applications}.
\newblock Mathematical Surveys and Monographs. Amer. Math.
Soc., Providence, Rhode Island (1988)

\bibitem{gamasa} 
A. G. Garc{\'i}a, F. Marcell\'{a}n, and L. Salto.
\newblock A distributional study of discrete classical
orthogonal polynomials.
\newblock {\em J. Comput. Appl. Math.} {\bf 57}, 147--162
(1995)

\bibitem{Hah} 
W. Hahn.
 \newblock Uber die Jacobischen polynome und zwei verwandte
Polynomklassen.
\newblock {\em Math. Zeit.} {\bf 39}, 634--638
(1935)

\bibitem{kost} 
R. Koekoek and R. F. Swarttouw.
\newblock {\em The Askey-scheme of hypergeometric orthogonal
polynomials and its q-analogue}, volume 98-17.
\newblock Reports of the Faculty of Technical Mathematics and
Informatics, Delft, The Netherlands (1998)

\bibitem{koo1} 
T. H. Koornwinder.
\newblock {\em Orthogonal polynomials in connection with
quantum groups}.
\newblock P. Nevai (Ed.), NATO ASI Series C 294. Kluwer Acad.
Publ., Dordrecht, The Netherlands,  257--292 (1990)

\bibitem{koo2} 
T. H. Koornwinder.
\newblock{\em Compact quantum groups and $q$-special
functions}.
\newblock V. Baldoni and M.A. Picardello (Eds.), Pitman
Research Notes in Mathematics Series. Longman Scientific
\& Technical, New York,  257--292 (1994)

\bibitem{koo4} 
T. H. Koornwinder.
\newblock The structure relation for Askey-Wilson
polynomials.
\newblock {\em J. Comput. Appl. Math.} {\bf 27}, 214--226
(2007)

\bibitem{mal} 
A. A. Malashin.
\newblock $q$-analog of Racah polynomials on the lattice
$x(s)=[s]q[s+1]_q$ and its connections with 6j-symbols
for the $su_q(2)$ and $su_q(1,1)$, quantum algebras.
\newblock Master's thesis, Moscow State University,
M.V. Lomonosov. In Russian  (1992)

\bibitem{mabrpe} 
F. Marcell\'{a}n, A. Branquinho, and J. C. Petronilho.
\newblock Classical orthogonal polynomials: a functional
approach.
\newblock {\em Acta Appl. Math.}  {\bf 34} , 283--303 (1994)

\bibitem{mealma} 
J. C. Medem, R. \'{A}lvarez-Nodarse, and F. Marcell\'{a}n.
\newblock On the $q$-polynomials: a distributional study.
\newblock {\em J. Comput. Appl. Math.}{\bf 135}, 157--196
(2001)

\bibitem{nag} 
Sh. M. Nagiyev.
\newblock Difference Schr\"{o}dinger equation and
$q$-oscillator model.
\newblock {\em Theoret. and Math. Phys.} {\bf 102},
180--187 (1995)

\bibitem{nisuuv} 
A. F. Nikiforov, S. K. Suslov, and V. B. Uvarov.
\newblock {\em Classical Orthogonal Polynomials of a Discrete
Variable}.
  Springer Series in Computational Physics.
\newblock Springer-Verlag, Berlin (1991)

\bibitem{niuv1} 
A. F. Nikiforov and V. B. Uvarov.
\newblock{\em The Special Functions of Mathematical Physics}.
\newblock Birkh\"{a}user Verlag, Basel (1988)

\bibitem{routh} 
E. Routh.
\newblock On some properties of certain solutions of a
differential equation of the second order.
\newblock {\em Proc. London Math. Soc.}{\bf 16}, 245--261
(1884)

\bibitem{rud} 
W. Rudin.
\newblock {\em Real and Complex Analysis}.
\newblock McGraw-Hill, New York (1974)

\bibitem{vikl} 
N. Ja. Vilenkin and A. U. Klimyk.
\newblock {\em Representations of Lie Groups and Special
Functions}, volume I, II, III.
\newblock Kluwer Academic Publishers, Dordrecht,
The Netherlands (1992)
\end{thebibliography}
\end{document}